\documentclass{amsart}

\theoremstyle{plain}

\newtheorem{thm}{Theorem}[section]

\newtheorem{lem}[thm]{Lemma}

\theoremstyle{definition}

\newtheorem{exmp}{Example}[section]

\theoremstyle{remark}

\newtheorem{rem}{Remark}[section]

\newcommand{\Rb}{\textrm{\bf R}\xspace}                             
\newcommand{\pscm}[1]{\ensuremath{\mathcal R^+(#1)}}        
\newcommand{\pscmt}[1]{\ensuremath{\mathcal R^+_0(#1)}}     
\newcommand{\bdr}[1]{\ensuremath{\partial #1}}              
\newcommand{\gra}{\ensuremath{\alpha}\xspace}

\newcommand{\grG}{\ensuremath{\Gamma}\xspace}

\newcommand{\gre}{\ensuremath{\epsilon}\xspace}
\newcommand{\grO}{\ensuremath{\Omega}\xspace}
\newcommand{\gro}{\ensuremath{\omega}\xspace}

\newcommand{\grf}{\ensuremath{\theta}\xspace}

\newcommand{\grt}{\ensuremath{\tau}\xspace}

\newcommand{\grk}{\ensuremath{\kappa}\xspace}

\newcommand{\grl}{\ensuremath{\lambda}\xspace}

\newcommand{\grd}{\ensuremath{\delta}\xspace}

\newcommand{\grx}{\ensuremath{\xi}\xspace}
\newcommand{\grPh}{\ensuremath{\Phi}\xspace}

\newcommand{\grps}{\ensuremath{\psi}\xspace}

\usepackage{latexsym,amsmath}
\usepackage{amssymb,xspace,amscd}

\begin{document}


\title[A quasifibration of psc metrics]{A quasifibration of spaces of positive scalar curvature metrics}
\author{Vladislav Chernysh}
\address{ Department of Mathematics\\ University of Notre
Dame\\Notre Dame, IN 46556}%
\email{vchernys@nd.edu}%

\begin{abstract} In this paper we show that for Riemannian
manifolds with boundary the natural restriction map is a
quasifibration between spaces of metrics of positive scalar
curvature. We apply this result to study homotopy properties of
spaces of such metrics on manifolds with boundary.
\end{abstract}

\maketitle


\section{Introduction}\label{intro}

The purpose of this note is to establish the fact that a natural
restriction map in Riemannian geometry is a quasifibration of
metrics of positive scalar curvature and discuss some applications
to the study of spaces of positive scalar curvature metrics.

If $M$ is an open manifold, then, by a result of Gromov
\cite{Gro69}, there always exists on $M$ a metric of positive
sectional curvature. However, such a metric, in general, will not
be complete. So, when studying metrics of positive scalar
curvature on manifolds with boundary, it is necessary to impose
some sort of  boundary condition. It is natural to require that a
metric is a product near the boundary.

Let $M$ be a manifold with boundary $\partial M$. We fix a collar
$c\colon \partial M\times (-1,0]\rightarrow M$ and define the
space $\pscm{M}$ of metrics of positive scalar curvature on $M$
that restrict to a product metric near the boundary with respect
to $c$. We take the usual Fr\'echet topology on this space. This
topology is defined by the collection of $C^k$-norms $\|.\|_k$ on
the space of all Riemannian metrics $\mathcal R (M^n)$ with
respect to some reference metric $h$: $\|g\|_k=\max_{i\le k}
\sup_{M^n}|\nabla^i g|$. The topology does not depend on the
choice of the metric $h$.

By $\pscmt{\partial M}$ we denote the image of $\pscm{M}$ under
the restriction map \[\rho\colon \pscm{M} \rightarrow
\pscm{\partial M},\]where $\rho(g):=g|_{\partial M}$. We assume
that $\pscm{M}$ is non empty, which, of course, implies that
$\pscmt{\partial M}$ is non empty.

The definition of a quasifibration is due to Dold and Thom
\cite{DT58}. A surjective map $p\colon E \rightarrow B$ is a
quasifibration if for all $x\in B$, all $y\in p^{-1}(x)$, and all
$i\ge 0$ we have $\pi_i(E,p^{-1}(x),y)\approx \pi_i(B,x)$.

\begin{thm}\label{main} The map $\rho\colon \pscm{M} \rightarrow \pscmt{\partial
M}$ is a quasifibration.
\end{thm}

To show that a map between topological spaces $f\colon
X\rightarrow Y$ (we take $Y$ to be path connected) is a
quasifibration it suffices to show that its fiber $f^{-1}(y_0)$ is
homotopy equivalent to its homotopy fiber under the canonical
inclusion map. The homotopy fiber of $f$ is defined by replacing
$f$ by the Serre path fibration $\hat{f}\colon \hat{X}\rightarrow
Y$ and taking the fiber $\grO_{y_0}:=\hat{f}^{-1}(y_0)$, some
$y_0\in Y$. Then $\grO_{y_0}=\{(x,\gro)\}$ where a path
$\gro\colon [0,1]\rightarrow Y$ is such that $\gro(0)=f(x)$ and
$\gro(1)=y_0$. Any two $\grO_{y_0}$ and $\grO_{y_1}$ are homotopy
equivalent. It follows that the homotopy type of the homotopy
fiber is well defined and we denote it as $\grO$.

The idea of the proof is to introduce an intermediate space
$\grO^s$, which is defined by taking smooth paths $\gro$ in the
definition of $\grO$. Then one can show that $\grO^s$ is homotopy
equivalent both to the fiber of $\rho$ and to the homotopy fiber
of $\rho$, see Lemmas~\ref{sdrtheorem} and~\ref{mil}.

One of important geometric implications of Theorem~\ref{main} is
Theorem~\ref{cont}.

A Hausdorff space $X$ is a topological manifold in the sense of
Palais \cite{Pal66} if there exists an open covering
$\{O_{\gra}\}$ of $X$ and a family of maps $\{\grf_{\gra}\colon
O_{\gra}\rightarrow V_{\gra}\}$, where each $V_{\gra}$ is locally
convex topological vector space and each $\grf_{\gra}$ is a
homeomorphism of $O_{\gra}$ onto either an open subset of
$V_{\gra}$ or an open subset of a half space of $V_{\gra}$.

\begin{thm}\label{cont} Let $A$ be a contractible subset of
$\pscmt{\partial M}$. Suppose that $A$ is a metrizable topological
manifold. Then for any point $a\in A$ the inclusion $i\colon
\rho^{-1}(a)\rightarrow \rho^{-1}(A)$ is a homotopy equivalence.
\end{thm}

\begin{rem} The conclusion of the above Theorem is also true under
the assumption that $\rho^{-1}(A)$ is an ANR or, more generally,
is dominated by a $CW$-complex. The author does not know whether
these properties follow from $A$ merely being contractible.
\end{rem}

In particular, if $h_0$ and $h_1$ are in the same path connected
component of $\pscmt{\partial M}$, then the spaces of positive
scalar curvature metrics that near the boundary restrict to a
product with correspondingly $h_0$ and $h_1$ are homotopy
equivalent.

Another geometric consequence of Theorem~\ref{main} is an
extension of results in \cite{Che04a} to manifolds with boundary.

Let $N^{n-k}\subset M^n$, $k\ge 3$, be a submanifold of $M^n$. We
assume that there exists a tubular neighborhood $\grt\colon
N\times D^k\rightarrow M$, such that the restriction $\grt\colon
\partial N\times D^k \rightarrow \partial M$ is a tubular
neighborhood of $\partial N$ in $\partial M$. We also assume that
the collar $c$ is compatible with $N$ in the sense that its
restriction to $\partial N\times (-1,0]$ is a collar for $\partial
N$.

We fix  a metric $g_N$ on $N$, and a torpedo metric $g_0$ on $D^k$
(a torpedo metric in the disc $D^k$ is an $O(k)$-symmetric,
positive scalar curvature metric, which is equal to a $k$-sphere
metric near the center of the disc and is a product with a
$(k-1)$-sphere metric near the boundary of the disc), such that
the metric $g_N+g_0$ has positive scalar curvature on $N\times
D^k$. Here the fixed metric $g_N$ can be any metric subject to the
only requirement that it is a product near the boundary $\partial
N$.

Let $h_0\in\pscmt{\partial M}$. Since codimension of $\partial N$
in $\partial M$ is greater than $2$, from~\cite{Che04a} we may
assume that $\grt^*(h_0)=g_N|_{\partial N}+g_0$. We define
\begin{equation*}\label{Ro}
(\rho^{-1}(h_0))_0:=\{g\in\rho^{-1}(h_0)|\grt^*(g)=g_N+g_0\}.
\end{equation*}

\begin{thm}\label{neat} Suppose that $\pscm{M}$ is not empty. Then
the inclusion map \[ i\colon (\rho^{-1}(h_0))_0\rightarrow
\rho^{-1}(h_0)\] is a homotopy equivalence.
\end{thm}

\begin{exmp}
Let $M^n$ be a manifold with a $k$handle $D^{n-k}\times D^k$, such
that $k\ge 3$ and $\pscm{M^n}$ is nonempty. Let $g_1$ be a metric
(which is a product near the boundary) on $D^{n-k}$ and $g_0$ is a
torpedo metric on $D^k$, such that $g:=g_1+g_0$ has positive
scalar curvature. Then for any $h_0\in\pscmt{\partial M}$ the
space $\rho^{-1}(h_0)$ is homotopy equivalent to the subspace of
$\rho^{-1}(\hat{h}_0)$ consisting of metrics that restrict to the
metric $g$ on the handle. Here the metric $\hat{h}_0$ is obtained
by deforming $h_0$ to be equal to $g_1|_{S^{n-k-1}}+ g_0$ on
$S^{n-k-1}\times D^k\subset\partial M$, see \cite{Che04a} for
details.
\end{exmp}

\section{Proofs of Theorems}

Given a smooth path of metrics $\gra\colon I\rightarrow \pscm{X}$
on a closed smooth manifold $X$, we would like to put a positive
scalar curvature metric on $X\times\Rb$. However, in general, the
scalar curvature of the obvious metric $g(x,t)=\gra(t)(x)+dt^2$
will not be positive.

We fix a smooth function $F\colon\Rb\rightarrow [0,1]$ such that
$0\le F'<2$, $F(t)=0,\, t\in(-\infty,\gre]$, $F(t)=1,\,
t\in[1-\gre,\infty)$, for some $0<\gre<1/4$, and for a positive
number $\grt$ we define a function $F_{\grt}\colon\Rb\rightarrow
\Rb$, by $F_{\grt}(t)=F(\frac {t}{\grt})$.

Let
\[
g^{\gra}_{\grt}(x,t):=\gra(F_{\grt}(t))(x)+dt^2.
\]

Define
\begin{eqnarray}
S'(\gra)&:=&\inf_{t>0}\{g^{\gra}_{\grt}\textrm{ is a psc metric
for all }\grt\ge t\},\nonumber\\
S(\gra)&:=&\max(S'(\gra),1).\label{S}
\end{eqnarray}

The function $S'(\gra)$ is clearly upper semi-continuous. By the
Lemma~\ref{Flemma} below, it is also lower semi-continuous. It
follows that $S$ is continuous and defines a metric on $X\times
\Rb$ by the formula
\begin{equation}\label{product}g^{\gra}_{t_0}(x,t):=\gra(F_{t_0}(t))(x)+dt^2,
\end{equation}
where $t_0:=S(\gra)$. By the same Lemma, this metric has positive
scalar curvature and is a Riemannian product near $X\times 0$ and
$X\times t_0$.

\begin{lem}\label{Flemma}

Let $\gra\colon[0,1]\rightarrow \pscm{X^n}$ be a
$C^{\infty}$-family of positive scalar curvature metrics on a
compact closed manifold $X^n$, then
\begin{itemize}
\item[(i)] $\exists\,\grl>0$ such that
$g^{\grl}(t):=\gra(F_{\grl}(t))+dt^2 \in \pscm{X^n\times[0,\grl]}$
and $g^{\grl}$ is a product metric near the boundary $(X\times
0)\cup(X\times 1)$ of $X\times [0,\grl]$;%

\item[(ii)] if $t_0=S(\gra)$ is a positive number, then
$\forall\, n=1,2,3,\dots\ \exists\, t_n>0,\,x_n\in X^n,\,\grt_n\in[0,t_n]$ such that, %
$t_0-\frac 1n<t_n<t_0 $ and the scalar curvature of
$\gra(F_{t_n}(t))+dt^2$ at $(x_n,\grt_n)$ is negative.

\end{itemize}

\end{lem}

\begin{proof} (i) Denote $g^{\grl}(x,t):=\gra(F_{\grl}(t))(x)+dt^2$. Let $(x_0,
\grt_0)$ be a point in $X^n\times[0,\grl]$. Take normal
coordinates for $\gra(F_{\grl}(\grt_0))$ at a point $x_0\in X^n$.
In these coordinates, we get $g^{\grl}_{ij}(x_0,
\grt_0)=\grd_{ij},\, \grG^k_{ij}=0$ for $1\le i,j,k \le n$. Recall
that
\[ \grG^k_{ij}=\frac 12 g^{kl}(\partial_i g_{jl} + \partial_j g_{il} -
\partial_l g_{ij}).\]
\noindent Since $g^{\grl}_{i,n+1} \equiv 0$ for $1\le i \le n$, we
get in our normal coordinates at $(x_0, \grt_0)$
\begin{eqnarray*}
\grG^k_{ij}(x_0,\grt_0)&=& 0 \quad \textrm{ for } 1\le i,j,k \le n\\
\grG^{n+1}_{n+1, i}&\equiv& 0 \quad \textrm{ for } 1\le i\le n+1\\
\grG^i_{n+1, j} (x_0, \grt_0) &=& \frac 12 \partial_{n+1}
g_{ij}(x_0,
\grt_0)\quad \textrm{ for } 1\le i,j \le n\\
\grG^{n+1}_{ij}(x_0, \grt_0) &=& -\frac 12 \partial_{n+1}
g_{ij}(x_0, \grt_0)\quad \textrm{ for } 1\le i,j \le n.
\end{eqnarray*}
\noindent And the equation for sectional curvature are
\begin{eqnarray*}
R^s_{ijk}&=&\partial_j \grG^s_{ik}-\partial_i
\grG^s_{jk}+\grG^l_{ik}\grG^s_{jl}-\grG^l_{jk}\grG^s_{il}\\
R_{ijks}&=&R^l_{ijk}g_{ls}\\
K_{ij}&=&(\partial_i,\partial_j, \partial_i,\partial_j)=R_{ijij}
\end{eqnarray*}
\noindent From the Gauss equation for curvature we get for $1\le
i,j \le n$
\[ K_{ij}=\overline{K}_{ij} + (b_{ii}b_{jj} - b^2_{ij}).\]
\noindent The remaining sectional curvatures
\begin{eqnarray*}
{K}_{i,n+1}&=& R^{n+1}_{i,n+1,i}\\
&=&\partial_{n+1}\grG^{n+1}_{ii}-\partial_i\grG^{n+1}_{n+1,i}+%
\grG^l_{ii}\grG^{n+1}_{n+1,l}-\grG^l_{n+1,i}\grG^{n+1}_{il}\\
&=&\partial_{n+1}\grG^{n+1}_{ii}-\grG^l_{n+1,i}\grG^{n+1}_{il}\\
&=& -\frac 12 \partial^2_{n+1}g_{ii}-(\frac 12
\partial_{n+1}g_{li})(-\frac 12\partial_{n+1}g_{il})\\
&=& -\frac 12 \partial^2_{n+1} g_{ii} + \frac 14
\sum_{l=1}^{n+1}(\partial_{n+1} g_{il})^2.
\end{eqnarray*}
\noindent Then the scalar curvature at a point $(x_0,\grt_0)$ is
given by the formula
\begin{eqnarray*}
\grk=\grk_X &+& \sum_{i,j=1}^n(b^2_{ij}-b_{ii}b_{jj})\\
&-&\sum_{i=1}^n\partial^2_{n+1}g_{ii} + \frac 12
\sum_{i,j=1}^n(\partial_{n+1}g_{ij})^2.
\end{eqnarray*}
\noindent Now, we have that $b_{ij}=\grG^{n+1}_{ij}=-\frac 12
\partial_{n+1}g_{ij}(x_0,\grt_0)$ and the formulas for the
derivatives
\begin{eqnarray*}
\partial_{n+1}g_{ij}(x_0,\grt_0)&=& \frac
1{t_0}F'(\grt_0)\gra'\left(F\left(\frac{\grt_0}{t_0}\right)\right)_{ij}(x_0)\\%
\partial^2_{n+1}g_{ij}(x_0,\grt_0) &=& \frac1{t_0^2}F''(\grt_0)\gra'%
\left(F\left(\frac{\grt_0}{t_0}\right)\right)_{ij}(x_0)\\%
&+&\frac1{t_0^2}(F'(\grt_0))^2\gra''\left(F\left(\frac{\grt_0}{t_0}\right)\right)_{ij}(x_0)
\end{eqnarray*}
\noindent The scalar curvature for the product may now be
expressed as
\begin{eqnarray*}
\grk=\grk_X &+& \frac 14 \frac 1{t^2_0}
\sum^n_{i,j=1}\left(F'\left(\frac{\grt_0}{t_0}\right)\right)^2%
\left((\gra_{ij}')^2 - \gra_{ii}'\gra_{jj}'\right)(x_0)\\
&+& \frac
1{t^2_0}\sum^n_{i=1}\left(F''\left(\frac{\grt_0}{t_0}\right)\gra_{ii}'(x_0)+
\left(F'\left(\frac{\grt_0}{t_0}\right)\right)^2\gra_{ii}''(x_0)\right)\\
&+& \frac 12 \frac 1{t^2_0}
\sum^n_{i,j=1}\left(F'\left(\frac{\grt_0}{t_0}\right)\right)^2\left(\gra_{ij}'(x_0)\right)^2
\end{eqnarray*}
To finish the proof, notice that $\grk_X$ is positive for all
$(x,\grt)\in X^n\times[0,t_0]$.\par%
\medskip\noindent
(ii) Suppose $t_0>0$ and let $(x_0, \grt_0)$ be a point in
$X^n\times[0,t_0]$ where the scalar curvature is not positive.
Such a point always exists since the $\grk>0$ is an open condition
and if the scalar curvature is everywhere positive we can find
$t_1<t_0$ such that the metric corresponding to $t_1$ will have
positive scalar curvature. Now , freeze the values of $t_0,\,
\grt_0,$ and $x_0$ which are in the arguments for the functions
$F,\, \gra$ and their derivatives, and regard the resulting
function as a function of the inverse of $t_0$. In the light of
the argument above it has a positive derivative at $t_0$, and its
value at $t_0$ is less or equal than $0$. So in an arbitrary
neighborhood \emph{on the left} from $t_0$ we can find a value
$t_n$ such that our function will be strictly negative at the
point $t_n$. Now, ``unfreezing" only $\grt_0$ we can find a number
$\grt_n$ such that $\frac{\grt_n}{t_n}=\frac{\grt_0}{t_0}$. The
point $(x_0, \grt_n)$ is the one that we were seeking.\end{proof}

We fix a metric $h_0\in\pscmt{\partial M}$ and consider the
homotopy fiber $\grO$ of $\rho$ at $h_0$, $\grO=\{(g,\gro)|\
\gro(0)=\rho(g),\, \gro(1)=h_0\}$. The topology on this fiber is
the usual compact open topology. The smooth homotopy fiber
$\grO^s$ is defined analogously by taking $\gro$ to be a smooth
path. We take the Fr\'echet topology on the smooth fiber.

There is a natural embedding $i$ of the fiber $\rho^{-1}(h_0)$
into $\grO^s$, $i(g)=(g,*)$, where $*$ is the constant path
$*(t)=h_0$.

\begin{lem}\label{sdrtheorem} The map $i\colon \rho^{-1}(h_0)\rightarrow\grO^s$
is a homotopy equivalence.
\end{lem}

\begin{proof} Let $(g,\gro)$ be a point in $\grO^s$ and $V_0$ be a constant
outward normal vector field on $\partial M$ of the unit length. We
take a smooth cutoff function $\grps$ on $\Rb$ with
$\grps(-3/4)=0$, $ \grps(-1/4)=1$, and define a vector field on $M$
by setting
\begin{eqnarray*}
V(x)&=&\begin{cases} \grps(t)S(\gro)V_0 &\quad
x=c(a,t)\\
0&\quad  \textrm{otherwise}
\end{cases},
\end{eqnarray*}
where $S$ is defined by the formula~\ref{S}. Extend this vector
field to $M\cup(\bdr{M}\times[0,\infty))$ as a constant vector
field $S(\gro)V_0$ on $\partial M\times[0,\infty)$ and denote by
$\grPh^S_1$ the diffeomorphism determined by the flow of this
vector field at $t=1$. Then
$\grPh^S_1(M)=M\cup(\bdr{M}\times[0,S])$. Define
\begin{eqnarray*}
g^{\gro}&=&\begin{cases} g &\quad \textrm{ on }M\\
g^{\gro}_{t_0} &\quad \textrm{ on }\bdr{M}\times [0,S]\\
h_0 &\quad \textrm{ on } \partial M\times [S,\infty)
\end{cases},
\end{eqnarray*}
\noindent where $g^{\gro}_{t_0}$ is given by the
formula~\ref{product}.

Now, we define an inverse map $r\colon \grO^s\rightarrow
\rho^{-1}(h_0)$ as
\[r(g,\gro):= \left(\grPh_1^{S(\gro)}\right)^{*}(g^{\gro}).
\] Here we take the restriction of the pullback metric to $M$.

For $u\in[0,1]$ we define a path
$\gro_u(\grt):=\gro((1-u)\grt+u)$. Let $g^{\gro}_{uF}$ be the
metric that is defined exactly as $g^{\gro}$ by taking the
function $uF$ instead of $F$. The homotopy $H\colon
\grO^s\times[0,1]\rightarrow \grO^s$ of $i\circ r$ to the identity
map is given by
\begin{eqnarray*}
H((g,\gro),u)=\begin{cases}
\left(\left(\grPh_1^{2uS(\gro)}\right)^{*}(g^{\gro}_{0F}),
\gro_{0}\right)& 0\le u\le 1/2,\\
\left(\left(\grPh_1^{S(\gro)}\right)^{*}(g^{\gro}_{(2u-1)F}),
\gro_{2u-1}\right)& 1/2 \le u\le 1.
\end{cases}
\end{eqnarray*}
When $u$ is equal to $0$, the map $H(\cdot,0)$ is the identity map
on $\grO^s$. When $u=1$, the map $H(\cdot,1)$ is equal to $i\circ
r$.\end{proof}

\begin{lem}\label{mil} The inclusion map $i\colon \grO^s\rightarrow \grO$ is
a homotopy equivalence.
\end{lem}

\begin{proof}The proof is completely analogous to the proof of Theorem~$17.1$
in~\cite{Mil63}. Since $\grO$ is an open subset of a locally
convex topological vector space, we can cover it with convex open
sets. Then take $\grO_k$, the space of all paths $\gro$ such that
$\gro\left([(j-1)/2^k,j/2^k]\right)$ is contained in some element
of the covering. The space $\grO$ is a homotopy direct limit of
$\grO_k$ and the space $\grO^s$ is a homotopy direct limit of
$\grO^s_k:=i^{-1}(\grO_k)$. By Milnor's argument, the map
\[i|_{\grO_k^s}\colon \grO^s_k\rightarrow \grO_k\] is a homotopy
equivalence. Here, the inverse map is defined by taking a path
$\gro\in\grO_k$ and assigning to it a piece-wise linear path that
coincides with $\gro$ at points $j/2^k$. Then we smooth the
resulting path by pre-composing with a smooth function that maps
$j/2^k$ to $j/2^k$ and all of whose derivatives vanish at points
$j/2^k$. This finishes the proof.\end{proof}

\begin{proof}[Proof of Theorem~\ref{cont}] If $p\colon
E\rightarrow B$ is a quasifibration over a contractible space $B$
then for any point $b\in B$ the inclusion of the fiber
$p^{-1}(b)\rightarrow E$ induces a weak homotopy equivalence. From
Palais~\cite{Pal66} we know that $\rho^{-1}(a)$ and $\rho^{-1}(A)$
are both dominated by $CW$-complexes. For such dominated spaces a
weak homotopy equivalence is, in fact, a homotopy equivalence by a
theorem of J.~H.~C.~Whitehead. From Theorem~\ref{main} it follows
that the inclusion map $i$ is a homotopy equivalence.\end{proof}

\begin{proof}[Proof of Theorem~\ref{neat}] As in the proof of
Theorem~\ref{cont}, it suffices to show that $i$ is a weak
homotopy equivalence. In \cite{Che04a} a method for deforming
compact families of metrics of positive scalar curvature was
developed, which allowed to prove the weak homotopy equivalence in
the case of closed manifolds $M$ and $N$. This deformation can be
readily adapted to manifolds with boundary and has an important
property. Namely, it preserves the product structure with respect
to the fixed tubular map $\grt$, i.e. if $\grt^*(g)=g_N+g_0$, then
for the deformation metrics $g(t)$, $t\in[0,1]$ we have
$\grt^*(g(t))=g_N+g_0(t)$ and $g(t)$ is constant outside of the
tubular neighborhood of $N$. The problem is that, in general,
$g_0(t)$ is not equal to $g_0$, so this deformation takes us
outside the fiber $\rho^{-1}(h_0)$.

The solution is to introduce a subspace $A\subset \pscmt{\partial
M}$ consisting of metrics that are equal to $h_0$ outside
$\grt(\partial N\times D^k)$ and equal to $g_N|_{\partial N}+g_w$
on $\grt(\partial N\times D^k)$. Here, $g_w$ is a warped metric in
the disc, i.e. $g_w=g(t)^2dt^2 + f(t)^2d\grx^2$, where $d\grx^2$
is the standard metric of the $(k-1)$-sphere of radius $1$, $g$ is
a smooth even function, and $f$ is a smooth odd function. Note
that $g_0\in A$, i.e. a torpedo metric is a warped metric. Then,
from the construction of the deformation, we have that $g_w(t)$ is
a warped metric for all $t\in [0,1]$. This allows us to conclude a
weak homotopy equivalence (and, therefore, a homotopy equivalence)
between $(\rho^{-1}(h_0))_0$ and $\rho^{-1}(A)$.

From \cite{Che04a} it follows that the inclusion map
$h_0\rightarrow A$ is a weak deformation retraction, cf.
Theorem~$4.1$ in \cite{Che04a}. Now, the proof follows from
Theorem~\ref{cont}.\end{proof}

\bibliographystyle{amsalpha}
\providecommand{\bysame}{\leavevmode\hbox
to3em{\hrulefill}\thinspace}
\providecommand{\MR}{\relax\ifhmode\unskip\space\fi MR }
\providecommand{\MRhref}[2]{%
  \href{http://www.ams.org/mathscinet-getitem?mr=#1}{#2}
} \providecommand{\href}[2]{#2}

\end{document}